\documentclass[graybox,footinfo]{svmult}

\smartqed
\usepackage{mathptmx}       % selects Times Roman as basic font
\usepackage{helvet}         % selects Helvetica as sans-serif font
\usepackage{courier}        % selects Courier as typewriter font
\usepackage{type1cm}        % activate if the above 3 fonts are
\usepackage{graphicx}       % standard LaTeX graphics tool
\usepackage{array,colortbl}
\usepackage{amsmath,amsfonts,amssymb,bm} % no amsthm, Springer defines Theorem, Lemma, etc themselves

\newcommand{\Z}{\mathbb{Z}} % integers
\newcommand{\R}{\mathbb{R}} % reals
\newcommand{\N}{\mathbb{N}} % natural numbers {1, 2, ...}
\newcommand{\rd}{\,\mathrm{d}} % differential symbol for use in integrals
\usepackage{microtype} % good font tricks

\usepackage[colorlinks=true,linkcolor=black,citecolor=black,urlcolor=black]{hyperref}
\usepackage{bookmark}
\makeatletter % to avoid hyperref warnings:
  \providecommand*{\toclevel@author}{999}
  \providecommand*{\toclevel@title}{0}
\makeatother

\begin{document}

%%%%%%%%%% OWN COMMANDS %%%%%%%%%%	

\newcommand{\novakbi}{{\boldsymbol{i}}}

\newcommand{\novakwidebar}[1]{\mbox{\kern1.5pt\hbox{\vbox{\hrule height 0.6pt \kern0.35ex
        \hbox{\kern-0.15em \ensuremath{#1 }\kern0.0em}}}}\kern-0.1pt}

\newcommand\novakCA{\emph{Constructive Approximation }}
\newcommand\novakFCM{\emph{Foundations of Computational Mathematics }}
\newcommand\novakJAT{\emph{Journal of Approximation Theory }}
\newcommand\novakJC{\emph{Journal of Complexity }} 
\newcommand\novakMC{\emph{Mathematics of Computation  }}
\newcommand\novakNM{\emph{Numerische Mathematik }}
\newcommand\novakPAMS{\emph{Proceedings of the American Mathematical Society }}
\newcommand\novakUSSR{\emph{USSR Computational Mathematics and Mathematical Physics }}

\title*{Some Results on the Complexity of Numerical Integration}
% if the original contribution title is too long you may use the following
% short title
% \titlerunning{   }
\author{Erich Novak}
\institute{Erich Novak
	\at Mathematisches Institut, University Jena, 
	Ernst-Abbe-Platz 2, D-07743 Jena, Germany\\
 \email{erich.novak@uni-jena.de}}
%The post address should be identical to that one of Robert J. Kunsch
\maketitle

\abstract{
We present 
some results on the complexity of numerical integration.  
We start with the seminal paper of Bakhvalov (1959) 
and end with new results 
on the curse of dimensionality and on the complexity of 
oscillatory integrals.\\
This survey paper consists of four parts: 
\begin{enumerate} 
	\item Classical results till 1971
	\item Randomized algorithms
	\item Tensor product problems, tractability and weighted norms
	\item Some recent results: $C^k$ functions and oscillatory integrals
\end{enumerate}
}

\section{Classical Results till 1971}  

I start with a warning: 
We do \emph{not} discuss the complexity of 
path integration and 
infinite-dimensional integration
on $\R^\N$ or other domains although 
there are exciting new results in that area,
see~\cite{novakBG14,novakCDGR07,novakDG14,novakDG14a,novakDKS13,novakGn12,
	novakG13,novakGMR13,novakHMNR09,novakKSWW09,novakNHMR11,novakNW10,novakPW11,
	novakWa14,novakWW96}. 
For parametric integrals see~\cite{novakDH13,novakDH14},
for quantum computers, see~\cite{novakHe01a,novakHe03a,novakN01,novakTW02}. 

We mainly study the problem of numerical integration, i.e., of approximating 
the integral
\begin{equation}\label{novakeq01}  
S_d(f) = \int_{D_d} f(x) \rd x
\end{equation}
over an open subset $D_d\subset \R^d$ of 
Lebesgue measure $\lambda^d(D_d)=1$ for integrable functions
$f\colon D_d\to\R$. 
The main interest is on the behavior of the minimal number of function values
that are needed in the worst case setting
to achieve an error at most $\varepsilon>0$. 
Note that classical examples of domains $D_d$ are the unit cube $[0,1]^d$ and 
the normalized Euclidean ball (with volume 1),
which are closed. However, we work with their interiors 
for definiteness of certain derivatives. 

We state our problem. Let $F_d$  be a class of %  continuous 
integrable functions $f\colon D_d\to\R$. For~\mbox{$f\in F_d$},
we approximate the integral $S_d(f)$, see~\eqref{novakeq01}, by algorithms
of the form 
$$ 
A_{n}(f)=\phi_{n}(f(x_1),f(x_2),\dots,f(x_n)),  
$$ 
where $x_j\in D_d$ can be chosen adaptively and $\phi_{n}\colon \R^n\to
\R$ is an arbitrary mapping. Adaption means that the selection of $x_j$
may depend on the already computed values $f(x_1),f(x_2),\dots,f(x_{j-1})$.
We define $N\colon  F_d \to \R^n$ by 
$N(f) = (f(x_1), \dots , f(x_n))$. 
The (worst case) error of the algorithm $A_{n}$ is defined by 
$$
e(A_{n})=\sup_{f\in F_d}|S_d(f)-A_{n}(f)|, 
$$
the optimal error bounds are given by
$$
e(n, F_d) = \inf_{A_n} \,  e(A_n) .
$$ 
The information complexity $n(\varepsilon,F_d)$ is 
the minimal number of function values 
which is needed to guarantee that the error is
at most $\varepsilon$, i.e., 
$$
n(\varepsilon,F_d)=\min\{n \mid \exists\ A_{n}\ \mbox{such that}\ 
e(A_{n})\le\varepsilon\}.
$$
We minimize $n$ over all choices of adaptive sample points $x_j$ and
mappings $\phi_{n}$. 

In this paper we give an overview on some of the basic results that 
are known about the numbers $e(n, F_d)$ and $n(\varepsilon, F_d)$. 
Hence we concentrate on complexity issues
and leave aside other important questions such as implementation 
issues. 

It was proved by Smolyak and Bakhvalov 
that as long as the class $F_d$
is convex and balanced 
we may restrict the minimization of $e(A_n)$ by
considering only nonadaptive 
choices of $x_j$ and  linear mappings $\phi_{n}$, 
i.e.,  it is enough to consider $A_n$ 
of the form
\begin{equation}     \label{novakeq03} 
A_n(f) = \sum_{i=1}^n a_i f(x_i) .
\end{equation}  

\setcounter{theorem}{-1} 
\begin{theorem}[Bakhvalov~\cite{novakBa71}]    \label{novakT0} 
Assume that the class $F_d$
is convex and balanced.  
Then 
\begin{equation} %R let it look nice
	e(n, F_d) = %\operatorname*{\vphantom{p}inf}_{x_1,\ldots,x_n}
			\inf_{x_1,\ldots,x_n}
				\sup_{\substack{f \in F_d \\ N(f) = 0}}
					S_d(f)
\end{equation} 
and for the infimum in the definition 
of $e(n, F_d)$ it is enough to consider 
linear and nonadaptive 
algorithms $A_n$ of the form \eqref{novakeq03}.
\end{theorem} 

In this paper we only
consider convex and balanced $F_d$ 
and then 
we can use the last formula for $e(n,F_d)$.

\setcounter{remark}{-1} 
\begin{remark} 
a) For a proof of Theorem~\ref{novakT0} see, 
for example, \cite[Theorem~4.7]{novakNW08}. This result is not really 
about complexity (hence it got its number), but it helps to 
prove complexity results. 

b) A~linear algorithm~$A_n$ is called a quasi Monte Carlo (QMC) algorithm 
if~\mbox{$a_i = 1/n$} for all~$i$ and is called a positive 
quadrature formula if~\mbox{$a_i > 0$} for all~$i$. 
In general it may happen that optimal quadrature formulas 
have some negative weights and, in addition, we cannot say much about the 
position of good points~$x_i$. 

c) More on the optimality of linear algorithms 
and on the power of adaption can be found in
\cite{novakCW04,novakNo96,novakNW08,novakTWW88,novakTW80}. 
There are important classes of functions that are \emph{not} 
balanced and convex, and where Theorem~\ref{novakT0} can not be applied, 
see also 
\cite{novakCDHHZ,novakPW09}.
\qed
\end{remark} 

The \emph{optimal order of convergence} plays an important role 
in numerical analysis.
We start with a classical result of Bakhvalov (1959) for the class 
$$
F^k_d = \{ f\colon  [0,1]^d \to \R \mid 
\Vert D^\alpha f \Vert_\infty \le 1,  \ 
|\alpha| \le k \},
$$
where $k \in \N$ and $|\alpha| = \sum_{i=1}^d \alpha_i $ for 
$\alpha \in \N_0^d$ and $D^\alpha f$ denotes the respective 
partial derivative. 
For two sequences $a_n$ and $b_n$ of positive numbers we write 
$a_n \asymp b_n$ if there are positive numbers $c$ and $C$ such 
that 
$c < a_n/b_n < C $ for all $n \in \N$. 

\begin{theorem}[Bakhvalov~\cite{novakBa59}]    \label{novakT1} 
\begin{equation} 
e(n, F_d^k) \asymp n^{-k/d}  .
\end{equation} 
\end{theorem} 

\begin{remark} %R m^d instead of k^d
a) For such a complexity result one needs to prove 
an upper bound (for a particular algorithm) 
and a lower bound (for all algorithms). 
For the upper bound one can use tensor product methods 
based on a regular grid, i.e., one can use the $n=m^d$ points $x_i$ 
with coordinates from the set $\{ 1/(2m), 3/(2m), \dots , (2m-1)/(2m) \}$. 

The lower bound can be proved with the technique 
of  ``bump functions'': One can construct $2n$ functions $f_1, \dots , f_{2n}$ 
with disjoint supports such that all $2^{2n}$ functions 
of the form 
$
\sum_{i=1}^{2n} \delta_i f_i 
$
are contained in $F^k_d$, where $\delta_i = \pm 1$ 
and $S_d(f_i) \ge c_{d,k} \, n^{-k/d-1}$.
Since an algorithm $A_n$ can only compute $n$ function values,
there are two functions $f^+ = \sum_{i=1}^{2n} f_i$ and 
$f^{-} = f^+ - 2   \sum_{k=1}^n f_{i_k} $ such that 
$f^+, f^- \in F^k_d$ and $A_n(f^+) = A_n(f^-)$ 
but $|S_d(f^+) - S_d(f^-)| \ge 2 n c_{d,k} n^{-k/d-1}$.
Hence the error of $A_n$ must be at least $c_{d,k} n^{-k/d}$. 
For the details see, for example, \cite{novakNo88}. 

b) Observe that we can not conclude much on 
$n(\varepsilon,F_d^k)$ if $\varepsilon$ is fixed and $d$ is large, since 
Theorem~\ref{novakT1} contains \emph{hidden factors} that depend on 
$k$ and $d$.  Actually the lower bound is of the form
$$
e(n, F_d^k) \ge c_{d,k} n^{-k/d} ,
$$
where the $c_{d,k}$ decrease with $d \to \infty$ and tend to zero. 

c) The proof of the upper bound (using tensor product algorithms) 
is easy since we assumed that the domain is $D_d = [0,1]^d$. 
The  optimal order of convergence is known for much more 
general spaces (such as Besov and Triebel-Lizorkin spaces) 
and arbitrary bounded Lipschitz domains, see~\cite{novakNT06,novakTr10,novakVy06}. 
Then the proof of the upper bounds is more difficult, however. 

d) Integration on fractals was recently studied 
by Dereich and M\"uller-Gron\-bach~\cite{novakDMG14}. 
These authors also obtain an optimal order of convergence 
$n^{-k/\alpha}$. The definition of $S_d$ 
must be modified and $\alpha$ coincides, under suitable conditions, 
with the Hausdorff dimension of the fractal. 
\qed
\end{remark} 

By the \emph{curse of dimensionality} we mean that 
$n(\varepsilon,F_d)$ is exponentially large in $d$. 
That is, there are positive numbers $c$, $\varepsilon_0$ and $\gamma$ such that
\begin{equation}   
n(\varepsilon,F_d) \ge c \, (1+\gamma)^d  \quad
\mbox{for all} \quad \varepsilon \le \varepsilon_0  \quad \mbox{and infinitely many} 
\quad  d\in \N. 
\end{equation}

If, on the other hand, 
$n(\varepsilon,F_d)$ is bounded by a polynomial in $d$ and $\varepsilon^{-1}$ 
then we say that the problem is \emph{polynomially tractable}. 
If 
$n(\varepsilon,F_d)$ is bounded by a polynomial in $\varepsilon^{-1}$ alone, i.e., 
$n(\varepsilon,F_d) \le C \varepsilon^{-\alpha}$ for $\varepsilon < 1$, 
then we say that the problem is \emph{strongly polynomially tractable}.

{}From the proof of Theorem~\ref{novakT1} 
we can not conclude whether 
the curse of dimensionality holds for the classes $F_d^k$ or not;
see Theorem~\ref{novakT11}. 
Possibly Maung Zho Newn and Sharygin~\cite{novakNS71} were the first 
who published (in 1971) 
a complexity result
for arbitrary $d$ with explicit constants and so proved the curse 
of dimensionality for Lipschitz functions. 

\begin{theorem}[Maung Zho Newn and Sharygin~\cite{novakNS71}]  \label{novakT2} 
Consider the class 
$$
F_d = \{ f \colon  [0,1]^d \to \R \mid 
|f(x)-f(y)| \le \max_i |x_i - y_i| \} .
$$
Then 
$$
e (n, F_d) = \frac{d}{2d+2} \cdot n^{-1/d} 
$$
for $n= m^d$ with $m \in \N$. 
\end{theorem}

\begin{remark} %R m^d instead of k^d
One can show that for $n=m^d$ the regular grid 
(points $x_i$ 
with coordinates from the set $\{ 1/(2m), 3/(2m), \dots , (2m-1)/(2m) \}$)
and the midpoint rule
$A_n(f) = n^{-1} \sum_{i=1}^n f(x_i)$ 
are optimal. 
See also \cite{novakBa76,novakBa77,novakCh95,novakSu79} for this result and for generalizations 
to similar function spaces.
\qed 
\end{remark}

\section{Randomized Algorithms}  

The integration problem is difficult 
for all deterministic algorithms if the classes $F_d$ 
of inputs are too large, see Theorem~\ref{novakT2}. 
One may hope that randomized algorithms make this problem much
easier. 

Randomized algorithms 
can be formalized in various ways leading to slightly different models. 
We do not explain the technical details and only give
a reason why it makes sense to study different models for upper and
lower bounds, respectively; see \cite{novakNW08} for more details. 

\begin{itemize}

\item
Assume that we want to construct and to analyze concrete algorithms that
yield upper bounds for the (total) complexity of given problems
including the arithmetic cost and the cost of generating random numbers.
Then it is reasonable to consider 
a rather restrictive model of computation where,
for example, only the standard arithmetic operations are allowed.
One may also restrict the use of random numbers and study so-called
\emph{restricted Monte Carlo methods}, where only random bits are allowed;
see \cite{novakHNP04}. 

\item
For the proof of lower bounds we take the opposite view
and allow  general randomized mappings
and a very general kind of randomness.
This makes the lower bounds stronger.

\end{itemize}

It turns out that the results are often very robust 
with respect to changes of
the computational model. For the purpose of this paper, 
it might be enough that a randomized algorithm $A$
is a random variable $(A^\omega)_{\omega \in \Omega}$ 
with a random element $\omega$ where, 
for each fixed $\omega$, the algorithm $A^\omega$ is a 
(deterministic) algorithm as before.
We denote by $\mu$ the distribution of the $\omega$. 
In addition one needs rather weak measurability assumptions, 
see also the textbook \cite{novakMNR09}. 
Let $\bar n(f,\omega)$ be the number of function values 
used for fixed $\omega$ and $f$. 

The number 
$$
\tilde n(A) = \sup_{f \in F} \int_\Omega \bar  n(f, \omega) \rd\mu (\omega)
$$
is called the \emph{cardinality} of the randomized algorithm $A$ and
$$
e^{\rm ran} (A) = \sup_{f \in F} \left( \int_\Omega^*   
\Vert S(f) - \phi_\omega (N_\omega (f))
\Vert^2 \rd\mu (\omega) \right)^{1/2}
$$
is the \emph{error} of $A$.
By $\int^*$ we denote the upper integral. For $n \in \N$,  define
$$
e^{\rm ran} (n, F_d) = \inf \{ e^{\rm ran} (A) \, : \, \tilde n(A) \le n \} .
$$

If $A\colon  F \to G$ is a (measurable) deterministic algorithm
then $A$ can also be treated as a randomized algorithm
with respect to a Dirac (atomic)  measure $\mu$.
In this sense we can say that deterministic algorithms are special
randomized algorithms. Hence the inequality
\begin{equation}
e^{\rm ran} (n, F_d) \le e (n, F_d)
\end{equation}
is trivial. 

The number $e^{\rm ran} (0, F_d)$ is called the 
\emph{initial error in the randomized setting}. 
For $n=0$, we do not sample $f$, and $A^\omega(f)$ is independent of
$f$, but may depend on~$\omega$. 
It is easy to check that for a linear $S$ 
and a balanced and convex set $F$, the best we can do is 
to take $A^\omega = 0$ and then
$$
e^{\rm ran}(0,F_d)=e(0,F_d).
$$
This means that for linear problems the initial errors are the same 
in the worst case and randomized setting.

The main advantage of randomized algorithms is that the curse 
of dimensionality is not present even for certain 
large classes of functions. 
With the standard Monte Carlo method we obtain 
$$
e^{\rm ran}(n,F_d) \le  \frac{1}{\sqrt{n}} , 
$$
when $F_d$ is the unit ball of $L_p([0,1]^d)$ and 
$2 \le p \le \infty$. 
Math\'e~\cite{novakMa95}  proved that this is almost optimal 
and the optimal algorithm is
$$
A^\omega_n(f) = \frac{1}{n+\sqrt{n}} \sum_{i=1}^n f(X_i) 
$$
with i.i.d. random variables $X_i$ that are uniformly distributed 
on $[0,1]^d$. It also follows that
$$
e^{\rm ran}(n,F_d) =    \frac{1}{1+ \sqrt{n}} , 
$$
when $F_d$ is the unit ball of $L_p([0,1]^d)$ and 
$2 \le p \le \infty$. 
In the case $1 \le p <2$ one can only achieve 
the rate $n^{-1+1/p}$, for a discussion see~\cite{novakHN02}. 

Bakhvalov~\cite{novakBa59} found the optimal order of convergence 
already in 1959 for the class 
$$
F^k_d = \{ f\colon  [0,1]^d \to \R \mid 
\Vert D^\alpha f \Vert_\infty \le 1,  \ 
 |\alpha| \le k \},
$$
where $k \in \N$ and $|\alpha| = \sum_{i=1}^d \alpha_i $ for 
$\alpha \in \N_0^d$. 

\begin{theorem}[Bakhvalov~\cite{novakBa59}]   
\begin{equation} 
e^{\rm ran} (n, F_d^k) \asymp n^{-k/d-1/2}  .
\end{equation} 
\end{theorem} 

\begin{remark} 
A~proof of the \emph{upper bound} can be given with a technique 
that is often called \emph{separation of the main part} 
or also \emph{control variates}. 
For $n=2m$ use $m$~function values to construct a good $L_2$ approximation 
$f_m$ of $f \in F_d^k$ by a deterministic algorithm. 
The optimal order of convergence is 
$$
\Vert f - f_m \Vert_2 \asymp m^{-k/d} .
$$
Then use the unbiased estimator
$$
A_n^\omega (f) = S_d (f_m) + \frac{1}{m} \sum_{i=1}^m (f-f_m)(X_i) 
$$ 
with i.i.d. random variables $X_i$ that are uniformly distributed 
on $[0,1]^d$. 
See, for example, \cite{novakMNR09,novakNo88} for more details. 
We add in passing that the optimal order of convergence can be obtained for 
many function spaces (Besov spaces, Triebel-Lizorkin spaces) 
and for arbitrary bounded Lipschitz domains
$D_d \subset \R^d$; see \cite{novakNT06}, where the approximation problem is 
studied. 
To obtain an explicit randomized 
algorithm with the optimal rate of convergence 
one needs a random number generator for the set $D_d$. 
If it is not possible to obtain 
efficiently random samples from the uniform distribution 
on $D_d$ one can work with Markov chain Monte 
Carlo (MCMC) methods, see Theorem~\ref{novakT5}. 

All known proofs of \emph{lower bounds} use the idea of Bakhvalov
(also called Yao's Minimax Principle): study the average case setting with 
respect to 
a probability 
measure on~$F$ and use the theorem of Fubini.  
For details
see~\cite{novakHe92,novakHe94,novakHe96,novakMNR09,novakNo88,novakNW10}. 
\qed
\end{remark} 

We describe a problem that was studied by several colleagues 
and solved by Hinrichs~\cite{novakHin10} using deep results from 
functional analysis. 
Let $H(K_d)$ be a reproducing kernel Hilbert space of real functions   
defined on a Borel measurable set $D_d\subseteq\R^d$.    
Its reproducing kernel $K_d:D_d\times   
D_d\to\R$ is assumed to be integrable,     
$$   
C^{\rm init}_d:=\left(\int_{D_d}\int_{D_d}K_d(x,y)\,\rho_d(x)\,\rho_d(y)   
\rd x \rd y\right)^{1/2}<\infty.   
$$   
Here, $\rho_d$ is a probability density function on $D_d$. Without   
loss of generality we assume that $D_d$ and $\rho_d$ are chosen such   
that there is no subset of $D_d$ with positive measure such that   
all functions from $H(K_d)$ vanish on it.    
   
The inner product and the norm of $H(K_d)$ are denoted by   
$\langle \cdot,\cdot\rangle_{H(K_d)}$ and $\|\cdot\|_{H(K_d)}$.   
Consider multivariate integration    
$$   
S_d(f)=\int_{D_d}f(x)\,\rho_d(x) \rd x\ \ \   
\mbox{for all}\ \ \ f\in H(K_d), 
$$   
where it is assumed that 
$S_d : H(K_d) \to \R$ is continuous. 

We approximate $S_d(f)$ in the randomized setting using   
\emph{importance sampling}. That is, for a positive 
probability density function   
$\tau_d$ on $D_d$    
we choose $n$ random sample points $x_1,x_2,\dots,x_n$ which   
are independent and distributed according to $\tau_d$ and take the   
algorithm   
$$   
A_{n,d,\tau_d}(f)=   
\frac1n\,\sum_{j=1}^n\frac{f(x_j)\,\rho_d(x_j)}{\tau_d(x_j)}.   
$$   
The error of $A_{n,d,\tau_d}$ is then 
$$   
e^{\rm ran}(A_{n,d,\tau_d})=\sup_{\|f\|_{H(K_d)}\le1}   
\left(\mathbb{E}_{\tau_d}\left(S_d(f)-A_{n,d,\tau_d}(f)   
\right)^2\right)^{1/2},   
$$   
where the expectation is with respect to the random choice of the   
sample points $x_j$.      
   
For $n=0$ we formally take $A_{0,d,\tau_d}=0$ and then   
$$   
e^{\rm ran}(0, H(K_d)) =C^{\rm init}_d.   
$$   

\begin{theorem}[Hinrichs~\cite{novakHin10}]   
Assume additionally that $K_d(x,y)\ge0$ for all $x,y\in D_d$.    
Then there exists a positive density function $\tau_d$ such that   
$$   
e^{\rm ran}(A_{n,d,\tau_d})\le \left(\frac{\pi}{2}\right)^{1/2}   
\frac1{\sqrt{n}}\ e^{\rm ran}(0,H(K_d)).   
$$   
Hence, if we want to achieve $e^{\rm ran}(A_{n,d,\tau_d})\le
\varepsilon \,e^{\rm ran}(0, H(K_d))$ then it is enough to take   
$$   
n=\left\lceil\frac{\pi}2\,\left(\frac1{\varepsilon}\right)^2\right\rceil.   
$$   
\end{theorem}   
 
\begin{remark} 
In particular, such problems are strongly polynomially   
tractable (for the normalized error) if the 
reproducing kernels are pointwise nonnegative and integrable. 
In \cite{novakNW11} we prove that the exponent $2$ of $\varepsilon^{-1}$ is sharp    
for tensor product Hilbert spaces whose univariate    
reproducing kernel is \emph{decomposable} and   
univariate integration is not trivial for the two    
parts of the decomposition. More specifically we have   
$$   
n^{\rm ran}(\varepsilon, H(K_d) )\ge\left\lceil\frac18\,   
\left(\frac1{\varepsilon}\right)^2\right\rceil   
\ \ \ \mbox{for all}\ \ \ \varepsilon \in(0,1)\ \    
\mbox{and}\ \ d\ge\frac{2\,\ln\,\varepsilon^{-1}\,-\,\ln\,2}{\ln\,\alpha^{-1}},    
$$   
where $\alpha\in[1/2,1)$ depends on the particular space.    
   
We stress that these estimates hold independently of the smoothness   
of functions in a Hilbert space. Hence, even for spaces of very smooth   
functions the exponent   
of strong polynomial tractability is~$2$.   
\qed 
\end{remark} 

Sometimes one cannot sample easily from the ``target distribution'' 
$\pi$ if one wants to compute an integral 
\[
S(f) = \int_D f(x)\, \pi(\mathrm{d} x).
\]
Then 
Markov chain Monte Carlo
(MCMC) methods are a very versatile and widely 
used tool.

We use an average of a finite Markov chain sample 
as approximation of the mean, i.e., we approximate $S(f)$ by
\[
  S_{n,n_0}(f) = \frac{1}{n} \sum_{j=1}^n f(X_{j+n_0}),
\]
where $(X_i)_{n \in \N_0}$ is a Markov chain with stationary 
distribution $\pi$. 
The number $n$ determines the number of
function evaluations of $f$.
The number $n_0$ is the \emph{burn-in}
or \emph{warm up} time. Intuitively, 
it is the number of steps of the 
Markov chain to get close to the stationary distribution $\pi$.

We study the mean square error of $S_{n,n_0}$, 
given by
\[
e_\nu(S_{n,n_0},f) = \left( \mathbb{E}_{\nu,K} 
\vert S_{n,n_0}(f)-S(f) \vert \right)^{1/2},
\]
where $\nu$ and $K$ indicate the initial distribution and the
transition kernel of the chain; 
we work with the spaces $L_p = L_p(\pi)$. 
For the proof of the following error bound we refer to 
\cite[Theorem~3.34 and Theorem~3.41]{novakRu12}.

\begin{theorem}[Rudolf~\cite{novakRu12}]   \label{novakT5}
Let $(X_n)_{n\in\mathbb{N}}$ be a Markov chain with reversible 
transition kernel $K$, 
initial distribution $\nu$, and 
transition operator $P$. 
Further, let
\begin{equation*}
  \Lambda = \sup\{ \alpha \colon \alpha \in \operatorname{{s}pec}(P-S) \},
\end{equation*}
where $\operatorname{{s}pec}(P-S)$ denotes 
the spectrum of the operator $(P-S) \colon L_2 \to L_2$,  and
assume that $\Lambda<1$.
Then
\begin{equation} 
    \sup_{\left \Vert f \right \Vert_{p}\leq 1 } e_\nu(S_{n,n_0},f)^2  
\leq \frac{2}{n(1-\Lambda)} + \frac{2\, C_\nu \gamma^{n_0}}{n^2(1-\gamma)^2}
\end{equation}
holds for $p=2$ and for $p=4$ under the following conditions: 
\begin{itemize} 
\item for $p=2$, $\frac{d\nu}{d\pi}\in L_\infty$ and a transition kernel
    $K$ which is $L_1$-exponentially convergent with $(\gamma,M)$ where
$\gamma <1$, i.e., 
$$
\Vert P^n -S \Vert_{L_1 \to L_1} \le  M \gamma^n 
$$
for all $n \in \N$ and 
    $ 
    C_\nu = M \left \Vert \frac{d\nu}{d\pi}-1\right \Vert_{\infty};
    $
\item 
for $p=4$, 
$\frac{d\nu}{d\pi}\in L_2$ and $\gamma= \Vert P-S \Vert_{L_2 \to 
L_2} <1$ where
    $
    C_\nu = 64 \left \Vert \frac{d\nu}{d\pi}-1\right \Vert_{2}.
    $
\end{itemize}
\end{theorem} 

\begin{remark} 
Let us discuss the results. 
First observe that we assume that the so called spectral gap 
$1-\Lambda$ is positive; in general we only know that 
$|\Lambda | \le 1$. 
If the transition kernel is $L_1$-exponentially convergent, then we
have an explicit error bound for integrands $f\in L_2$ whenever the initial 
distribution has a density $\frac{d \nu}{ d \pi} \in L_\infty$.
However, in general it is difficult to provide 
explicit values $\gamma$ and $M$ such that
the transition kernel is $L_1$-exponentially 
convergent with $(\gamma,M)$. This motivates 
to consider transition kernels which satisfy a 
weaker convergence property, such as the existence
of an $L_2$-spectral gap, i.e., 
$\Vert P-S \Vert_{L_2 \to L_2} <1$.  
In this case we have an 
explicit error bound for integrands $f\in L_4$ 
whenever the initial distribution has a 
density $\frac{d \nu}{ d \pi} \in L_2$. Thus,
by assuming a weaker convergence property of 
the transition kernel we obtain a weaker result in the sense
that $f$ must be in $L_4$ rather than $L_2$. 
 
If we want to have an error of $\varepsilon \in(0,1)$ it is still not
clear how to choose $n$ and $n_0$ to minimize 
the total amount of steps $n+n_0$. 
How should we choose 
the burn-in $n_0$?
One can prove in this setting, see \cite{novakRu12}, 
that the choice 
$ n^* 
= \lceil\frac{\log C_\nu}{1-\gamma} \rceil$ 
is a reasonable and almost optimal choice for the burn-in.

More details can be found in \cite{novakNR14}. 
For a full discussion with all the proofs see \cite{novakRu12}. 
\qed
\end{remark} 

\section{Tensor Product Problems and Weights}  

We know from the work of Bakhvalov 
already done in 1959
that the optimal order of convergence is $n^{-k/d}$ for functions 
from the class $C^k([0,1]^d)$. To obtain 
an order of convergence of roughly $n^{-k}$ for every dimension $d$,
one needs stronger smoothness conditions. 
This is a major reason for the study of functions with \emph{bounded mixed 
derivatives}, or \emph{dominating mixed smoothness},   such as the classes
$$
W^{k, {\rm mix}}_p ([0,1]^d) = 
\{ f: [0,1]^d \to \R \mid 
\Vert D^\alpha f \Vert_p \le 1 \hbox{ for } 
\Vert \alpha \Vert_\infty \le k\} .
$$

Observe that functions from this class have, in particular, 
the high order derivative $D^{(k,k, \dots , k)}f \in L_p$ 
and one may hope that the curse of dimensionality 
can be avoided or at least moderated by this assumption. 
For $k=1$ these spaces are closely related to various 
notions of \emph{discrepancy}, 
see, for example, \cite{novakDKS13,novakDP08,novakLP14,novakNW10,novakTem03a}. 

The optimal order of convergence is known for all $k \in \N$ 
and $1 < p < \infty$ due to the work of 
Roth~\cite{novakRo54,novakRo80},
Frolov~\cite{novakFr76,novakFr80}, Bykovskii~\cite{novakBy85}, 
Temlyakov~\cite{novakTem90} and Skriganov~\cite{novakSk95}, 
see the survey
Temlyakov~\cite{novakTem03a}. 
The cases $p \in \{ 1, \infty \}$ are still unsolved. 
The case $p=1$ is strongly related to the star discrepancy, see 
also Theorem~\ref{novakT-star}. 

\begin{theorem}    \label{novakt6} 
Assume that $k \in \N$ and $1< p < \infty$. 
Then 
$$
e(n, W_p^{k, {\rm mix}} ([0,1]^d)  ) \asymp n^{-k} (\log n)^{(d-1)/2} . 
$$
\end{theorem} 

\begin{remark} 
The upper bound was proved by Frolov~\cite{novakFr76} for $p = 2$ 
and by Skriganov~\cite{novakSk95} for all $p>1$. 
The lower bound was proved by Roth~\cite{novakRo54} and Bykovskii \cite{novakBy85}  
for $p=2$ and by Temlyakov~\cite{novakTem90} for all $p < \infty$. 
Hence it took more than 30 years to prove Theorem~\ref{novakt6} completely.

For functions in 
$W_p^{k, {\rm mix}} ([0,1]^d)$ 
with compact support in $(0,1)^d$ 
one can take algorithms of the form 
$$
A_n(f) = \frac{|\operatorname{det} A |}{a^d} 
\sum_{m \in {\mathbb Z}^d} f \left( \frac{Am}{a} \right) ,
$$
where $A$ is a suitable matrix that does not depend on $k$ or $n$, 
and $a>0$. 
Of course the sum is finite since we use only the points 
$ \frac{Am}{a} $ in $(0,1)^d$. 

This algorithm is similar to a 
\emph{lattice rule} but is not quite a lattice rule
since the points do not build an integration lattice. 
The sum of the weights is roughly 1, but not quite. 
Therefore this algorithm is not really a 
quasi-Monte Carlo algorithm. 
The algorithm $A_n$ can be modified to obtain the optimal order 
of convergence for the whole space
$W_p^{k, {\rm mix}} ([0,1]^d)$. 
The modified algorithm uses different points $x_i$ 
but still positive weights $a_i$.  
For a tutorial on this algorithm see \cite{novakUl14}. 
Error bounds for Besov spaces are studied in 
\cite{novakD97}.  Triebel-Lizorkin spaces and the case of small smoothness 
are studied in \cite{novakUU15} and \cite{novakNUU15}. 
\qed

For the Besov-Nikolskii classes $S^{r}_{p,q} B (T^d)$ 
with $1 \le p,q \le \infty$ and $1/p < r < 2$, the optimal 
rate is 
$$
n^{-r} (\log n) ^{(d-1)(1-1/q)} 
$$
and can be obtained constructively with QMC algorithms, 
see~\cite{novakHMOU}. 
The lower bound was proved by Triebel~\cite{novakTr10}. 
\qed

The Frolov algorithm can be used as a building block for a 
randomized algorithm that is \emph{universal} in the sense that 
it has the optimal order of convergence 
(in the randomized setting as well as in the worst case setting) 
for many different function spaces, see \cite{novakKN15}. 
\qed 
\end{remark} 

A~famous algorithm for \emph{tensor product problems}
is the \emph{Smolyak algorithm}, also called 
\emph{sparse grids algorithm}. 
We can mention just a few papers and books 
that deal with this topic: 
The algorithm was invented by
Smolyak~\cite{novakSm63} and, independently, 
by several other colleagues and research groups. 
Several error bounds were proved by 
Temlyakov~\cite{novakTem87,novakT93}; 
explicit error bounds (without unknown constants) 
were obtained by 
Wasilkowski and Wo\'zniakowski~\cite{novakWW95,novakWW99}.
Novak and Ritter~\cite{novakNR96b,novakNR97,novakNR99} 
studied the particular Clenshaw-Curtis Smolyak algorithm.
A~survey is Bungartz and Griebel~\cite{novakBG04} and another 
one is~\cite[Chap.~15]{novakNW10}. 
For recent results on the order of convergence see 
Sickel and T.~Ullrich~\cite{novakSU06,novakSU11} and
Dinh D\~ung  and T.~Ullrich~\cite{novakDU2014}. 
The recent paper~\cite{novakHNU14} contains a tractability result 
for the Smolyak algorithm applied to very smooth functions. 
We display only one recent result on the Smolyak 
algorithm. 

\begin{theorem}[Sickel and T.~Ullrich~\cite{novakSU11}]  
For the classes 
$W_2^{k, {\rm mix}} ([0,1]^d)$ 
one can construct a 
Smolyak algorithm
with the order of the error 
\begin{equation}    \label{novakbound-sickel} 
n^{-k} (\log n)^{(d-1)(k+1/2)} . 
\end{equation} 
\end{theorem} 

\begin{remark} 
a) The bound \eqref{novakbound-sickel} is valid 
even for $L_2$ approximation instead of integration, 
but it is not known whether this upper bound is 
optimal for the approximation problem.
Using the technique of control variates one can obtain the order
$$ 
n^{-k-1/2} (\log n)^{(d-1)(k+1/2)} 
$$
for the integration problem in the randomized setting. 
This algorithm is not often used since it is not easy 
to implement and its arithmetic cost is rather high. 
In addition, the rate can be improved by 
the algorithm of \cite{novakKN15} to
$ 
n^{-k-1/2} (\log n)^{(d-1)/2}
$.

b) 
It is shown in Dinh D\~ung  and T.~Ullrich~\cite{novakDU2014}
that the order \eqref{novakbound-sickel} 
can not be improved when restricting to Smolyak grids.

c) 
We give a short description of the Clenshaw-Curtis Smolyak 
algorithm for the computation of 
integrals $\int_{[-1,1]^d} f(x) \rd x$  that 
often leads to ``almost optimal'' 
error bounds, see \cite{novakNR97}. 

We assume that for $d=1$ a sequence of formulas 
$$
U^i(f) = \sum_{j=1}^{m_i} a_j^i \,  f(x_j^i) 
$$
is given. In the case of numerical integration the $a_j^i$ 
are just numbers. 
The method $U^i$ uses $m_i$ function values and we assume 
that 
$U^{i+1}$ has smaller error than $U^i$ and $m_{i+1} > m_i$. 
Define then, for $d>1$, the tensor product formulas
$$
(U^{i_1} \otimes \cdots \otimes U^{i_d})(f) = 
\sum_{j_1=1}^{m_{i_1}} \!\cdots\!  \sum_{j_d=1}^{m_{i_d}} 
a_{j_1}^{i_1} \cdots  a_{j_d}^{i_d} \;
f(x_{j_1}^{i_1} , \dots  , x_{j_d}^{i_d} ) . 
$$
A~tensor product formula clearly needs
$$
m_{i_1} \cdot m_{i_2} \cdot \dots \cdot m_{i_d} 
$$
function values, sampled on a regular grid. 
The Smolyak formulas
$A(q,d)$ are clever linear combinations of
tensor product formulas such that

\begin{itemize} 

\item
only tensor products with a relatively small number of
knots are used;

\item 
the linear combination is chosen in such a way that
an interpolation property for $d=1$ is preserved for 
$d>1$.

\end{itemize} 

The Smolyak formulas are defined by 
$$
A(q,d) = 
\sum_{q-d+1 \le |\novakbi| \le q} (-1)^{q-|\novakbi|} \cdot 
\binom{d-1}{q-|\novakbi|} \cdot
(U^{i_1} \otimes \cdots \otimes U^{i_d}), 
$$
where $q \geq d$.
Specifically, we use,  
for $d>1$, the Smolyak construction and start, 
for $d=1$, with the classical 
Clenshaw-Curtis formula with 
$$
m_1=1 \quad \hbox{  and } \quad  m_i = 2^{i-1}+1 
\ \hbox{ for } \ i>1. 
$$
The Clenshaw-Curtis formulas 
$$
U^i(f) = \sum_{j=1}^{m_i} a_j^i \,  f(x_j^i)
$$
use the knots  
$$
x_j^i= -\cos \frac{\pi (j-1)}{m_i - 1},   \qquad j = 1, \dots , m_i 
$$
(and $x_1^1=0$). 
Hence we use nonequidistant knots. 
The weights $a_j^i$ are defined in such a way
that $U^i$ is exact for all (univariate) polynomials of 
degree at most $m_i$.  \qed 
\end{remark} 

It turns out that many tensor product problems are still 
intractable and suffer from the curse of dimensionality, 
for a rather exhaustive presentation see~\cite{novakNW08,novakNW10,novakNW12}. 
Sloan and Wo\'zniakowski~\cite{novakSW98} describe a very interesting idea 
that was further  developed  in hundreds of papers, the paper 
\cite{novakSW98} is most important and influential. 
We can describe here only the very beginnings of a long 
ongoing  story; we present just one example instead of the whole theory. 

The rough idea is that $f \colon [0,1]^d \to \R$ may depend on many variables,
$d$ is large, but some variables or groups of 
variables are more important than others. 
Consider, for $d=1$, the inner product
$$ %R with parentheses to make it look less confusing
\langle f,g \rangle_{1, \gamma} 
= \left(\int_0^1  f \rd x\right) \left(\int_0^1  g  \rd x\right) + \frac{1}{\gamma} 
\int_0^1 f'(x) \, g'(x) \rd x ,
$$
where $\gamma >0$. 
If $\gamma $ is small then $f$ must 
be ``almost constant'' if it has small norm. 
A~large $\gamma$ means that $f$ may have a large variation and still the norm 
is relatively small. 
Now we take tensor products of such spaces and weights 
$\gamma_1 \ge \gamma_2 \ge \dots$ 
and consider the complexity of the integration problem 
for the unit ball $F_d$ with respect to this weighted norm. 
The kernel $K$ of the tensor product space $H(K)$ is 
of the form
$$
K (x,y) = \prod_{i=1}^d K_{\gamma_i} (x_i, y_i) ,
$$
where $K_\gamma$ is the kernel of the respective space $H_\gamma$
of univariate functions. 

\begin{theorem}[Sloan and Wo\'zniakowski~\cite{novakSW98}]   \label{novakT8} 
Assume that 
$\sum_{i=1}^\infty  \gamma_i < \infty$. 
Then the problem 
is strongly polynomially tractable. 
\end{theorem} 

\begin{remark}
The paper \cite{novakSW98} contains also a lower bound 
which is valid for all quasi-Monte Carlo methods. 
The proof of the \emph{upper bound}  is very interesting and an 
excellent example for the 
\emph{probabilistic method}. 
Compute the mean of the quadratic 
worst case error of  QMC algorithms  over all 
$(x_1, \dots , x_n) \in [0,1]^{nd}$ and obtain
$$
\frac{1}{n} \left( \int_{[0,1]^d} K (x,x) \rd x
	- \int_{[0,1]^{2d}} K(x,y) \rd x \rd y \right) .
$$ 
This expectation is of the form $C_d \,  n^{-1}$ 
and the sequence $C_d$ is bounded if and only if 
$\sum \gamma_i < \infty$. 
The 
\emph{lower bound} in \cite{novakSW98} is based on the fact that the kernel 
$K$ is always non-negative; this leads to lower bounds for
QMC algorithms 
or, more generally, for algorithms with positive weights.   \qed 

As already indicated, 
Sloan and Wo\'zniakowski~\cite{novakSW98} was continued in many directions. 
Much more general weights and many different Hilbert spaces 
were studied. 
By the probabilistic method one only obtains the 
\emph{existence} of a good QMC algorithms but, in the meanwhile, 
there exist many results about the \emph{construction} 
of good algorithms. 
In this paper the focus is on the basic complexity results and therefore 
we simply list a few of the most relevant papers: 
\cite{novakBDLNP12, novakSC02, novakDLPW09, novakDSWW04, novakDSWW06,  
novakHSW04b,        
novakHW2000, 
novakHW01, 
novakKPW12, 
novakKPW14, 
novakK03,  
novakKWWat06,    
novakNC06a,       
novakNC06b,       
novakSKJ02,  
novakSR02,
novakSWW}. 
See also the books 
\cite{novakDP08,novakLP14,novakNi92,novakNW10} 
and the excellent survey paper
\cite{novakDKS13}.   \qed 
\end{remark}

In complexity theory we want to study 
\emph{optimal} algorithms and it is not clear whether 
QMC algorithms or quadrature formulas with positive 
coefficients $a_i$ are optimal. 
Observe that the Smolyak algorithm uses also negative 
$a_i$  and it is known that 
in certain cases positive quadrature formulas are far from 
optimal; for examples see
\cite{novakNSW97} or \cite[Sects.~10.6 and~11.3]{novakNW10}. 
Therefore it is not clear whether the conditions 
on the weights in Theorem~\ref{novakT8}  can be relaxed if we allow 
arbitrary algorithms. 
The next result shows that this is not the case.

\begin{theorem}[\cite{novakNW99}]   \label{novakT9}
The integration problem from Theorem~\ref{novakT8} 
is strongly polynomially tractable
if and only if \, 
$\sum_{i=1}^\infty  \gamma_i < \infty$. 
\end{theorem} 

\begin{remark}
Due to the known upper bound of Theorem~\ref{novakT8}, to prove Theorem~\ref{novakT9}  
it is enough to prove a \emph{lower} bound for arbitrary algorithms. 
This is done via the technique of 
\emph{decomposable kernels} that was developed in \cite{novakNW99},
see also \cite[Chap.~11]{novakNW10}. 

We do not describe this technique here and only remark 
that we need for this technique many non-zero 
functions $f_i$ in the Hilbert space $F_d$ with disjoint 
supports. 
Therefore this technique usually works for functions with 
finite smoothness, but not for analytic functions.
\qed
\end{remark} 

Tractability of integration can be proved for many 
weighted spaces
and one may ask whether there are also 
unweighted spaces where tractability holds as well. 
A~famous example for this are integration problems 
that are related to the \emph{star discrepancy}. 

For $x_1 , \dots , x_n \in [0,1]^d$ define the star discrepancy by 
$$ 
D^*_\infty (x_1, \dots , x_n) = \sup_{ t \in [0,1]^d}
\left| % \rule{0cm}{0.8cm}\right. 
t_1 \cdots t_d - \frac{1}{n} 
\sum_{i=1}^n 1_{[0,t)}(x_i)
% \left.\rule{0cm}{0.8cm}
\right| , 
$$ 
the respective QMC quadrature formula is 
$Q_n(f)= \frac{1}{n} \sum _{i=1}^n f(x_i)$. 

Consider the Sobolev space
$$
F_d = \{ f \in W_1^{1, {\rm mix}}  \mid \Vert f \Vert \le 1, \,  f(x)=0
\hbox{ if there exists an $i$ with } x_i=1 \}
$$
with the norm 
$$
\Vert f \Vert  := \left\Vert \frac{\partial^d f}{\partial x_1 \partial x_2 
\dots \partial x_d} \right\Vert_1 .
$$
Then the 
Hlawka-Zaremba-equality yields 
$$
D^*_\infty (x_1, \dots , x_n)
= \sup_{f \in F_d} |S_d(f)-Q_n(f) |,
$$
hence the star discrepancy is a worst case error bound for integration. 
We define 
$$
n(\varepsilon, F_d ) = \min \{ n \mid \exists  \, x_1, \dots , x_n \hbox{ with } 
D^*_\infty (x_1, \dots , x_n) \le \varepsilon \}.
$$
The following result shows that this integration 
problem is polynomially tractable 
and the complexity is linear in the dimension. 

\begin{theorem}[\cite{novakHNWW99}]    \label{novakT-star} 
\begin{equation}    \label{novakstar} 
n( \varepsilon , F_d) \le C \,  d  \,  \varepsilon^{-2} 
\end{equation} 
and 
$$
n(1/64, F_d) \ge 0.18 \, d . 
$$ 
\end{theorem}

\begin{remark} 
This result was modified and improved in various ways and we mention 
some important results. 
Hinrichs~\cite{novakHinr} proved the lower bound  
$$ 
n( \varepsilon , F_d ) \ge c \, d \,  \varepsilon^{-1} \quad \hbox{for} 
\quad \varepsilon \le \varepsilon_0. 
$$
Aistleitner~\cite{novakA11} proved that the constant $C$ in \eqref{novakstar} 
can be taken as 100. 
Aistleitner and Hofer~\cite{novakAH14} proved more on upper bounds. 
Already the proof in \cite{novakHNWW99} showed that an upper bound 
$D_\infty^* (x_1, \dots , x_n) \le C \,  \sqrt{\frac{d}{n}}$ 
holds with high probability if the points $x_1, \dots, x_n$ 
are taken independently and uniformly distributed. 
Doerr~\cite{novakDo14} proved the respective lower bound, hence
$$
\mathbb{E}  (D_\infty^* (x_1, \dots , x_n)) 
\asymp \sqrt{\frac{d}{n}} \quad \hbox{for} \quad n \ge d. 
$$
Since the upper bounds 
are proved with the probabilistic method, we only know the 
\emph{existence} of points with small star discrepancy. 
The existence results can be 
transformed into (more or less explicit)  constructions 
and the problem is, of course, 
to minimize the computing time as well as the 
discrepancy. 
One of the obstacles is that already the computation of 
the star discrepancy of given 
points $x_1, x_2, \dots , x_n$ is very difficult. 
We refer the reader 
to~\cite{novakD08,novakDP13,novakDP14,novakDOG06,novakDGKP07,novakDGW10,
	novakDGW14,novakGn12a,novakHin12}. \qed

Recently Dick~\cite{novakD14} proved a tractability result 
for another unweighted space that 
is defined via an $L_1$-norm and consists of periodic 
functions; we denote Fourier coefficients by 
$\tilde f(k)$, where $k \in \Z^d$. 
Let $0 < \alpha \le 1$ and $1 \le p \le \infty$
and 
$$ %R \sup_{x,h}
F_{\alpha, p, d} = \left\{ f: [0,1]^d \to \R \mid 
\sum_{k \in \Z^d} | \tilde f (k) | 
+ \sup_{x,h} \frac{|f(x+h) -f(x)|}{\Vert h \Vert^\alpha_p} \le 1 
\right\}.
$$
Dick proved the upper bound 
$$
e(n, F_{\alpha, p, d})  
\le \max \left( \frac{d-1}{\sqrt{n}}, \frac{d^{\alpha/p}}{n^\alpha} 
\right) 
$$
for any prime number $n$. 
Hence the complexity is at most quadratic in $d$. 

The proof is constructive, a suitable algorithm is the following.
Use points \mbox{$x_k = \left( \left\{ \frac{k^1}{n} \right\}  , 
\left\{  \frac{k^2}{n}  \right\}  , \dots ,
\left\{ \frac{k^d}{n} \right\} 
\right)$}, where  \mbox{$  k= 0, 1, \dots , n-1$},
and take the respective QMC algorithm. 
\qed
\end{remark} 

\section{Some Recent Results}  

We end this survey with two results that were still unpublished 
at the time of the conference, April 2014. 
First we return to the classes $C^k([0,1]^d)$, see Theorem~\ref{novakT1}. 
We want to be a little more general and consider the 
computation of 
\begin{equation}  
S_d(f) = \int_{D_d} f(x) \rd x  %  \quad \mbox{for} \quad  f\in F_d 
\end{equation} 
up to some error $\varepsilon >0$, where  
$D_d \subset \R^d$ has Lebesgue measure 1. 
The results hold for arbitrary sets $D_d$,  
the standard example of course is $D_d = [0,1]^d$.   
For convenience we 
consider functions 
$f \colon \R^d \to \R$. 
This makes the function class a bit smaller
and the result a bit stronger, 
since our emphasis is on lower bounds. 
 
It has \emph{not} been 
known if the curse of dimensionality is present for probably the most  
natural class which is the unit ball  
of $r$ times continuously differentiable functions,  
\[ 
F_d^k 
=\{f\in C^k( \R^d ) \mid \Vert D^\alpha f \Vert_\infty \le 1 \quad 
\mbox{for all} \quad  |\alpha|\le k\}, 
\] 
where $k \in \N$. 

\begin{theorem}[\cite{novakHNUW12}]   \label{novakT11} 
The curse of dimensionality holds for the classes  
$F_d^k$ with the  
\emph{super-exponential} lower bound 
\[ 
n(\varepsilon, F_d^k ) \ge c_k\,(1-\varepsilon) \, d^{\,d /(2k+3)} 
\quad 
\text{for all}  \ d\in\N  
\ \text{and} \ \varepsilon\in(0,1),  
\]  
where $c_k >0$ depends only on $k$. 
\end{theorem} 

\begin{remark} 
In \cite{novakHNUW12,novakHNUW13} we 
also prove that the curse of dimensionality holds 
for even smaller classes of functions~$F_d$ for which the norms of  
arbitrary directional derivatives are bounded proportionally to $1/\sqrt{d}$. 
 
We start with the fooling function 
\[ 
f_0 (x) = \min\left\{1, \frac{1}{\delta\sqrt{d}}\, 
\operatorname{dist}(x,\mathcal{P}_\delta)\right\}  
\quad \mbox{for all} \quad x\in\R^d,  
\] 
where  
\[ 
\mathcal{P}_\delta = \bigcup_{i=1}^n B_\delta^d(x_i) 
\] 
and $B_\delta^d(x_i)$ is the ball with center $x_i$ and radius  
$\delta\sqrt{d}$.  
The function $f_0$ is Lipschitz. By a suitable smoothing  
via convolution we construct a smooth fooling function  
$f_k \in F_d$ with $f_k|_{\mathcal{P}_0} = 0$.  

Important elements of the proof are volume estimates 
(in the spirit of Elekes~\cite{novakE86} 
and Dyer, F\"uredi and McDiarmid~\cite{novakDFM90}), since
we need that the volume of a neighborhood of the convex hull
of $n$ arbitrary points is exponentially small  in $d$.   \qed 

Also classes of $C^\infty$-functions were studied recently. 
We still do not know whether 
the integration problem suffers from the curse of 
dimensionality for the classes
$$
F_d = \{ f: [0,1]^d \to \R \mid \Vert D^\alpha f \Vert_\infty \le 1 
\hbox{ for all } \alpha \in \N_0^d \}, 
$$
this is Open Problem 2 from \cite{novakNW08}. 
We know from 
Vyb{\'\i}ral~\cite{novakVyb13}
and \cite{novakHNUW13} 
that the curse is present for somewhat larger spaces 
and that a weak tractability holds for smaller 
classes;
this can be proved with the Smolyak algorithm, 
see~\cite{novakHNU14}. 
\qed
\end{remark} 

We now consider 
univariate oscillatory integrals for       
the standard Sobolev spaces      
$H^s$ of periodic and non-periodic      
functions with an arbitrary integer $s\ge1$.      
We study the approximate computation of
Fourier coefficients
$$ 
I_k (f) = \int_0^1 f(x) \, \eul^{-2\pi\,\imag\,kx} \rd x,        
\qquad \imag=\sqrt{-1},        
$$ 
where $k \in \Z$ and $f \in H^s$.  

There are several recent papers         
about the approximate computation of highly oscillatory  
univariate integrals      
with the weight $\exp(2\pi\,\imag\,kx)$, where $x\in[0,1]$ and $k$  
is an integer (or $k \in \R$) which  is      
assumed to be large in the absolute sense, see      
Huybrechs and Olver~\cite{novakHO09} for a survey.      

We study the Sobolev space $H^s$ for a finite $s\in\N$, i.e.,         
\begin{equation}         
H^s = \{ f: [0,1] \to \mathbb{C} \mid         
f^{(s-1)} \hbox{ is abs. cont., } f^{(s)} \in L_2 \}         
\end{equation}        
with the inner product        
\begin{equation} 
\begin{split}      
\langle f,g \rangle_{s} \;&=\; \sum_{\ell=0}^{s-1} \int_0^1 f^{(\ell)}(x) \rd  x  \;      
        \int_0^1 \novakwidebar{g^{(\ell)}(x)} \rd x  \,+\, \int_0^1 f^{(s)}(x)\,      
\novakwidebar{ g^{(s)}(x)} \rd x  \\      
&=\; \sum_{\ell=0}^{s-1} \langle f^{(\ell)},1\rangle_0 \,      
\novakwidebar{\langle g^{(\ell)},1\rangle_0}       
                        \,+\, \langle f^{(s)}, g^{(s)} \rangle_0,      
\end{split}      
\end{equation}      
 where         
$\langle f,g \rangle_0 = \int_0^1 f(x)\, \novakwidebar{g(x)} \rd x$,      
and norm $\|f\|_{H^s}=\langle f,f \rangle_{s}^{1/2}$.     

For the periodic case,
an  algorithm that uses $n$ function values at equally      
spaced points is nearly optimal,  
and its worst case error is bounded by $C_s(n+|k|)^{-s}$       
with $C_s$ exponentially small in $s$.       
For the non-periodic case, we first compute      
successive derivatives up to order $s-1$ at the end-points      
$x=0$ and $x=1$. These derivatives values are used to periodize the      
function and this allows us to obtain similar error bounds like for the      
periodic case.
Asymptotically in $n$, the worst case error      
of the algorithm is of order $n^{-s}$ independently of $k$       
for both periodic and non-periodic cases.       

\begin{theorem}[\cite{novakNUW13}]   
Consider the integration problem $I_k$ defined over the        
space $H^s$ of non-periodic functions with $s \in \N$.        
Then       
$$        
\frac{c_s}{(n+|k|)^s} \;\le\; e(n,k,H^s) \;\le\;        
\left(\frac3{2\pi}\right)^s \frac{2}{(n+|k|-2s+1)^s} ,        
$$        
for all $k\in\Z$ and $n\ge 2s$.    
\end{theorem}         

\begin{remark} 
The minimal errors $e(n,k,H^s)$ for the non-periodic case have a       
peculiar property for $s\ge2$ and large $k$. Namely,        
for $n=0$ we obtain the initial error which is of order $|k|^{-1}$,       
whereas for $n\ge 2s$ it becomes of order $|k|^{-s}$. Hence, the       
dependence on $|k|^{-1}$ is short-lived and disappears quite quickly.       
For instance, take $s=2$. Then $e(n,k,H^s)$ is of order $|k|^{-1}$        
only for $n=0$ and maybe for $n=1,2,3$, and then 
becomes of order $|k|^{-2}$.        
\qed
\end{remark} 

%%%%%%%%%%%%%%%%%%%%%%%%%%%%%%%%%%%%%%%%%%%%%%%%%%%%%%%%%%%%%%%%%%%%%%%%%%%%%%%%

\begin{acknowledgement}
I thank the following colleagues and friends for valuable 
remarks: 
Michael Gnewuch,
Aicke Hinrichs, 
Robert Kunsch, 
Thomas M\"uller-Gron\-bach,
Daniel Rudolf, 
Tino Ullrich, 
and
Henryk Wo\'zniakowski. 
I also thank two referees for carefully reading my manuscript. 
\end{acknowledgement}

%%%%%%%%%%%%%%%%%%%%%%%%%%%%%%%%%%%%%%%%%%%%%%%%%%%%%%%%%%%%%%%%%%%%%%%%%%%%%%%%%%%%%%%%%%%
%%% The bibliography
%
% BibTeX users please use
%\bibliographystyle{spmpsci}
%\bibliography{mybibfile}

\begin{thebibliography}{99.}%

\bibitem{novakA11}
Ch.~Aistleitner.
\newblock Covering numbers, dyadic chaining and discrepancy.
\newblock \novakJC 27:531--540, 2011.

\bibitem{novakAH14}
Ch.~Aistleitner and M.~Hofer.
\newblock Probabilistic discrepancy bounds for Monte Carlo point sets.
\newblock \novakMC 83:1373--1381, 2014.

\bibitem{novakBa76}
V.~F.\ Babenko.
\newblock Asymptotically sharp bounds for the remainder for the best
quadrature formulas for several classes of functions.
%\newblock \emph{Mat. Zametki} 19(3):313--322, 1976. 
\newblock English translation: \emph{Mathematical Notes} 19(3):187--193, 1976. 

\bibitem{novakBa77}
V.~F.\ Babenko.
\newblock 
Exact asymptotics of the error of weighted cubature formulas optimal 
for certain classes of functions. 
%\newblock \emph{Mat. Zametki} 20(4):589--595, 1976. 
\newblock English translation: \emph{Mathematical Notes} 20(4):887--890, 1976. 

\bibitem{novakBa59}       
N.~S.\ Bakhvalov.
\newblock On the approximate calculation of multiple integrals.
\newblock {\em Vestnik MGU, Ser. Math. Mech. Astron. Phys. Chem.}        
4:3--18, 1959, in Russian.       
\newblock English translation: \novakJC 31:502--516, 2015. 

\bibitem{novakBa71}
N.~S.\ Bakhvalov.
\newblock On the optimality of linear methods for operator        
approximation in convex classes of functions.        
\newblock \novakUSSR  11:244--249, 1971.        
       
\bibitem{novakBDLNP12}
J.~Baldeaux, J.~Dick, G.~Leobacher, D.~Nuyens, and F.~Pillichshammer.
\newblock Efficient calculation of the worst-case error and (fast) 
component-by-component construction 
of higher order polynomial lattice rules.
\newblock \emph{Numerical Algorithms} 59:403--431, 2012. 

\bibitem{novakBG14}
J.~Baldeaux and M.~Gnewuch.
\newblock Optimal randomized multilevel algorithms for 
infinite-dimensional integration on function spaces with 
ANOVA-type decomposition.
\newblock \emph{SIAM Journal Numerical Analysis} 52:1128--1155, 2014. 

\bibitem{novakBG04}       
H.-J.\ Bungartz and M.~Griebel.
\newblock Sparse grids.
\newblock \emph{Acta Numerica} 13:147--269, 2004.        

\bibitem{novakBy85}
V.~A.\ Bykovskii.
\newblock On the correct order of the error of optimal cubature formulas 
in spaces with dominant derivative, and on quadratic deviations of grids.
\newblock \emph{Computing Center Far-Eastern Scientific 
Center, Akad. Sci. USSR, Vladivostok}, preprint, 1985.

\bibitem{novakSC02}       
W.~W.~L.\ Chen and M.~M.\ Skriganov.
\newblock Explicit constructions in the classical mean squares problem in       
irregularities of point distribution.
\newblock {\em Journal f\"ur  Reine und  Angewandte Mathematik (Crelle)} 545:67--95, 2002.       

\bibitem{novakCh95}
E.~V.\ Chernaya.
\newblock Asymptotically exact estimation of the error 
of weighted cubature formulas optimal 
in some classes of continuous functions. 
\newblock \emph{Ukrainian Mathematical Journal } 47(10):1606--1618, 1995.
       
\bibitem{novakCDHHZ}
N.~Clancy, Y.~Ding, C.~Hamilton, F.~J.\ Hickernell, and Y.~Zhang. 
\newblock The cost of deterministic, adaptive, automatic algorithms: 
Cones, not balls.
\newblock \novakJC  30:21--45, 2014. 

\bibitem{novakCDGR07}
J.~Creutzig, S.~Dereich, Th.~M\"uller-Gronbach, and K.~Ritter.
\newblock Infinite-dimensional quadrature and approximation of        
distributions.
\newblock \novakFCM 9:391--429, 2009.    
       
\bibitem{novakCW04}
J.~Creutzig and P.~Wojtaszczyk.
\newblock Linear vs.\ nonlinear algorithms for linear problems.
\newblock \novakJC 20:807--820, 2004.       

\bibitem{novakDH13}
T.~Daun and S.~Heinrich.
\newblock Complexity of Banach space valued and parametric integration.
\newblock In J.~Dick, F.~Y.\ Kuo, G.~W.\ Peters, and I.~H.\ Sloan, editors,
	{\em {M}onte {C}arlo and Quasi-{M}onte {C}arlo Methods 2012},
	pages 297--316. Springer-Verlag, 2013.

%R
\bibitem{novakDH14}
T.~Daun and S.~Heinrich.
\newblock Complexity of parametric integration in various smoothness classes.
\newblock Preprint, 2014.
%\newblock \novakJC 30:750--766, 2014.

%R
\bibitem{novakDMG14}
S.~Dereich and Th.~M\"uller-Gronbach.
\newblock Quadrature for self-affine distributions on $\R^d$.
\newblock To appear in \novakFCM. 

\bibitem{novakD08}
J.~Dick.
\newblock A note on the existence of sequences with small star discrepancy.
\newblock \novakJC 23:649--652, 2007.

\bibitem{novakD14}       
J.~Dick.
\newblock Numerical integration of H\"older continuous, absolutely convergent 
Fourier-, Fourier cosine-, and Walsh series.
\newblock \novakJAT 183:14--30, 2014. 

\bibitem{novakDG14}
J.~Dick and M.~Gnewuch.
\newblock Optimal randomized changing dimension algorithms 
for infinite-dimensional integration on function spaces with 
ANOVA-type decomposition.
\newblock \novakJAT 184:111--145, 2014.

\bibitem{novakDG14a}
J.~Dick and M.~Gnewuch.
\newblock Infinite-dimensional integration in weighted 
Hilbert spaces: anchored decompositions, optimal 
deterministic algorithms, 
and higher order convergence.
\newblock \novakFCM 14:1027--1077, 2014. 

\bibitem{novakDKS13}       
J.~Dick,  F.~Y.\ Kuo, and I.~H.\ Sloan.
\newblock High-dimensional integration: 
The quasi-Monte Carlo way.
\newblock \emph{Acta Numerica} 22:133--288, 2013. 

\bibitem{novakDLPW09}
J.~Dick, G.~Larcher, F.~Pillichshammer, and H.~Wo\'zniakowski.   
\newblock Exponential convergence and tractability of multivariate
integration for Korobov spaces.
\newblock \novakMC 80:905--930, 2011.
   
\bibitem{novakDP08}
J.~Dick and  F.~Pillichshammer.
\newblock {\em Digital Nets and Sequences: Discrepancy Theory and Quasi-Monte
  Carlo Integration}.       
\newblock Cambridge University Press, 2010. 

%R: The template contains only a proceedings as an example
\bibitem{novakDP13}       
J.~Dick and F.~Pillichshammer.
\newblock Discrepancy theory and quasi-Monte Carlo integration.
\newblock In W.~Chen, A.~Srivastav, and G.~Travaglini,
editors, \emph{Panorama in Discrepancy Theory},
pages 539--619. Lecture Notes in Mathematics 2107, Springer-Verlag, 2014.

\bibitem{novakDP14}       
J.~Dick and F.~Pillichshammer.
\newblock The weighted star discrepancy of Korobov's $p$-sets.
\newblock To appear in \novakPAMS.       

\bibitem{novakDSWW04}        
J.~Dick, I.~H.\ Sloan, X.~Wang, and H.~Wo\'zniakowski.
\newblock Liberating the weights.
\newblock \novakJC 20:593--623, 2004.

\bibitem{novakDSWW06}        
J.~Dick, I.~H.\ Sloan, X.~Wang, and H.~Wo\'zniakowski.       
\newblock Good lattice rules in weighted Korobov spaces with general weights.
\newblock \novakNM 103:63--97, 2006.        

\bibitem{novakDU2014} 
Dinh~D\~ung and T.~Ullrich.
\newblock Lower bounds for the integration error for multivariate functions with mixed 
smoothness and optimal Fibonacci cubature for functions on the square.
\newblock \emph{Mathematische Nachrichten} 288:743--762, 2015.  

\bibitem{novakDo14}
B.~Doerr.
\newblock A lower bound for the discrepancy of a random point set.
\newblock \novakJC  30:16--20, 2014.

\bibitem{novakDOG06}        
B.~Doerr and M.~Gnewuch.
\newblock Construction of low-discrepancy point sets of        
small size by bracketing covers and dependent randomized        
rounding.
\newblock In A.~Keller, S.~Heinrich, and H.~Niederreiter, editors,
\emph{Monte Carlo and Quasi-Monte Carlo Methods 2006},
pages 299--312.        
Springer-Verlag, 2008.        
       
\bibitem{novakDGKP07}       
B.~Doerr, M.~Gnewuch, P.~Kritzer, and F.~Pillichshammer.
\newblock Component-by-component construction of low-discrepancy        
point sets of small size.
\newblock \emph{Monte Carlo Methods and Applications } 14:129--149, 2008.        
  
\bibitem{novakDGW10}
B.~Doerr, M.~Gnewuch, and M.~Wahlstr\"om.
\newblock Algorithmic construction of low-discrepancy 
point sets via dependent randomized rounding.
\newblock \novakJC 26:490--507, 2010. 

%R no example in the template that totally fits to this type of reference.
\bibitem{novakDGW14}
C.~Doerr, M.~Gnewuch, and M.~Wahlstr\"om.
\newblock Calculation of discrepancy measures and applications.
\newblock %to appear in: 
In W.~W.~L.\ Chen, A.~Srivastav, and G.~Travaglini, editors, 
\emph{Panorama of Discrepancy Theory}, pages 621--678.
Lecture Notes in Mathematics 2107, Springer-Verlag, 2014.

\bibitem{novakD97}
V.~V.\ Dubinin.
\newblock Cubature formulas for Besov classes.
\newblock \emph{Izvestija Mathematics } 61(2):259--283, 1997. 

%R No template, shortened the reference to the standard format
\bibitem{novakDFM90}
M.~E.\ Dyer, Z.~F{\"u}redi, and C.~McDiarmid.
\newblock Random volumes in the {$n$}-cube.       
\newblock {DIMACS Series in Discrete Mathematics and
Theoretical Compututer Science} 1:33--38, 1990.

\bibitem{novakE86}       
G.~Elekes.
\newblock A geometric inequality and the complexity of computing volume.       
\newblock \emph{Discrete Computational Geometry } 1:289--292, 1986.

\bibitem{novakFr76}
K.~K.\ Frolov.
\newblock Upper bounds on the error of quadrature formulas
on classes of functions.
\newblock \emph{Doklady Akademy Nauk USSR } 231:818--821, 1976.       
\newblock English translation: Soviet Mathematics Doklady 17, 1665--1669, 1976. 

\bibitem{novakFr80}       
K.~K.\ Frolov.
\newblock Upper bounds on the discrepancy in $L_p$, $2 \le p <\infty$,      
\newblock \emph{Doklady Akademy Nauk USSR} 252:805--807, 1980. 
\newblock English translation: Soviet Mathematics Doklady 18(1):37--41, 
1977.
       
\bibitem{novakGn12}
M.~Gnewuch.
\newblock Infinite-dimensional integration on weighted Hilbert spaces.
\newblock \novakMC 81:2175--2205, 2012.

\bibitem{novakGn12a}
M.~Gnewuch.
\newblock Entropy, randomization, derandomization, and discrepancy.
\newblock In L.~Plaskota and H.~Wo\'zniakowski, editors, 
\emph{Monte Carlo and Quasi-Monte Carlo Methods 2010}, pages 43--78.
Springer-Verlag, 2012. 

\bibitem{novakG13}
M.~Gnewuch.
\newblock Lower error bounds for randomized multilevel and 
changing dimension algorithms.
\newblock In J.~Dick, F.~Y.\ Kuo, G.~W.\ Peters, and I.~H.\ Sloan, editors, {\em {M}onte {C}arlo
  and Quasi-{M}onte {C}arlo Methods 2012}, pages 399--415. Springer-Verlag, 2013. 

\bibitem{novakGMR13}
M.~Gnewuch, S.~Mayer, and K.~Ritter.
\newblock On weighted Hilbert spaces and integration of functions 
of infinitely many variables.
\newblock \novakJC 30:29--47, 2014.

\bibitem{novakHe92}       
S.~Heinrich.
\newblock Lower bounds for the complexity of Monte Carlo function        
approximation.
\newblock \novakJC 8:277--300, 1992.        

%R no template with similar reference
\bibitem{novakHe94}     
S.~Heinrich.
\newblock Random approximation in numerical analysis.
\newblock In K. D. Bierstedt et al., editors, \emph{Functional Analysis},        
pages 123--171. Dekker, 1994.

%R no tomplate with similar reference
\bibitem{novakHe96}       
S.~Heinrich.
\newblock Complexity of Monte Carlo algorithms.
\newblock In \emph{The Mathematics of Numerical Analysis}, pages 405--419.
Lectures in Applied Mathematics 32, AMS-SIAM Summer Seminar, Park City,        
American Mathematical Society, 1996.        

\bibitem{novakHe01a}       
S.~Heinrich.
\newblock Quantum Summation with an Application to Integration.       
\newblock \novakJC 18:1--50, 2001.        

\bibitem{novakHe03a}
S.~Heinrich.
\newblock Quantum integration in Sobolev spaces.
\newblock \novakJC 19:19--42, 2003.        

\bibitem{novakHN02}
S.~Heinrich and E.~Novak.
\newblock Optimal summation and integration by deterministic, 
randomized, and quantum algorithms.
In K.-T.\ Fang, F.~J.\ Hickernell, and H.~Niederreiter, editors,
{\em {M}onte {C}arlo
  and Quasi-{M}onte {C}arlo Methods 2000}, pages 50--62.
Springer-Verlag, 2002. 

\bibitem{novakHNP04}       
S.~Heinrich, E.~Novak, and  H.~Pfeiffer.
\newblock How many random bits do we need for Monte Carlo integration?       
\newblock In H. Niederreiter, editor,
\emph{Monte Carlo and Quasi-Monte Carlo Methods 2002}, pages 27--49.
Springer-Verlag, 2004.
       
\bibitem{novakHNWW99}       
S.~Heinrich, E.~Novak, G.~W.\ Wasilkowski, and H.~Wo\'zniakowski.
\newblock The inverse of the star-discrepancy depends linearly on the dimension.        
\newblock \emph{Acta Arithmetica} 96:279--302, 2001.        

%R K.~Ritter instead of Klaus~Ritter?
\bibitem{novakHMNR09}
F.~J.\ Hickernell, Th.~M\"uller-Gronbach, B.~Niu, and K.~Ritter.
\newblock Multi-level Monte Carlo algorithms for    
infinite-dimensional integration on $\R^\N$.
\newblock \novakJC 26:229--254, 2010.

\bibitem{novakHSW04b}        
F.~J.\ Hickernell, I.~H.\ Sloan, and G.~W.\ Wasilkowski.       
\newblock On strong tractability of weighted multivariate integration.
\newblock \novakMC 73:1903--1911, 2004.

\bibitem{novakHW2000}       
F.~J.\ Hickernell and H.~Wo{\'z}niakowski.
\newblock Integration and approximation in arbitrary dimension.
\newblock \emph{Advances of Computational Mathematics } 12:25--58, 2000.

\bibitem{novakHW01}       
F.~J.\ Hickernell and H.~Wo{\'z}niakowski.
\newblock Tractability of multivariate integration for periodic functions.
\newblock \novakJC 17:660--682, 2001.

\bibitem{novakHinr}       
A.~Hinrichs.
\newblock Covering numbers, Vapnik-Cervonenkis classes and bounds for the star       
discrepancy.
\newblock \novakJC 20:477--483, 2004.        

\bibitem{novakHin10}   
A.~Hinrichs.
\newblock Optimal importance sampling for the approximation of integrals.
\novakJC 26:125--134, 2010.

\bibitem{novakHin12}   
A.~Hinrichs.
\newblock Discrepancy, integration and tractability.
\newblock In J.~Dick, F.~Y.\ Kuo, G.~W.\ Peters, and I.~H.\ Sloan, editors, {\em {M}onte {C}arlo
  and Quasi-{M}onte {C}arlo Methods 2012}, pages 129--172. Springer-Verlag, 2013.

\bibitem{novakHMOU}
A.~Hinrichs, L.~Markhasin,~J. Oettershagen, and T.~Ullrich.
\newblock Optimal quasi-Monte Carlo rules on higher order digital nets for 
the numerical integration of multivariate periodic functions.
\newblock Submitted.

\bibitem{novakHNUW12}
A.~Hinrichs, E.~Novak, M.~Ullrich, and H.~Wo\'zniakowski.
\newblock The curse of dimensionality for numerical integration
of smooth functions.
\newblock \novakMC 83:2853--2863, 2014.  

\bibitem{novakHNUW13}
A.~Hinrichs, E.~Novak, M.~Ullrich, and H.~Wo\'zniakowski.
\newblock The curse of dimensionality for numerical integration
of smooth functions II.
\newblock \novakJC 30:117--143, 2014.

\bibitem{novakHNU14}
A.~Hinrichs, E.~Novak, and M.~Ullrich.
\newblock On weak tractability of the Clenshaw Curtis Smolyak 
algorithm.
\newblock \novakJAT 183:31--44, 2014. 

\bibitem{novakHO09}
D.~Huybrechs and S.~Olver.
\newblock Highly oscillatory quadrature.
\newblock \emph{London Mathematical Society Lecture Note Series} 366:25--50, 2009.

\bibitem{novakKN15}
D.~Krieg and E.~Novak.
\newblock A universal algorithm for multivariate integration. 
\newblock Manuscript, available at arXiv. 
  
\bibitem{novakKPW12}
P.~Kritzer, F.~Pillichshammer, and H.~Wo\'zniakowski.
\newblock Multivariate integration of infinitely many times differentiable
functions in weighted Korobov spaces. 
\newblock \novakMC 83:1189--1206, 2014.  

\bibitem{novakKPW14}
P.~Kritzer, F.~Pillichshammer, and H.~Wo\'zniakowski.
\newblock Tractability of multivariate analytic problems.
\newblock In \emph{Uniform distribution and quasi-Monte Carlo methods}, 
pages 147--170. 
De Gruyter, 2014. 

\bibitem{novakK03}       
F.~Y.\ Kuo.
\newblock Component-by-component constructions achieve the optimal rate        
of convergence for multivariate integration in weighted Korobov and        
Sobolev spaces.
\newblock \novakJC 19:301--320, 2003.

\bibitem{novakKSWW09}   
F.~Y.\ Kuo, I.~H.\ Sloan, G.~W.\ Wasilkowski, and H.~Wo\'zniakowski.
\newblock Liberating the dimension.
\newblock \novakJC 26:422--454, 2010.
   
\bibitem{novakKWWat06}
F.~Y.\ Kuo, G.~W.\ Wasilkowski, and B.~J.\ Waterhouse.
\newblock Randomly shifted lattice rules for unbounded integrands.
\newblock \novakJC 22:630--651, 2006.    
  
\bibitem{novakLP14}
G.~Leobacher and F.~Pillichshammer.
\newblock \emph{Introduction to Quasi-Monte Carlo Integration and Applications}.
\newblock Springer-Verlag, 2014. 
       
\bibitem{novakMa95}      
P.~Math\'e.
\newblock The optimal error of {M}onte {C}arlo integration.
\newblock \novakJC 11:394--415, 1995.       
       
\bibitem{novakMNR09}    
Th.~M\"uller-Gronbach, E.~Novak, and K.~Ritter.
\newblock \emph{Monte-Carlo-Algorithmen}.
\newblock Springer-Verlag, 2012.

%R no example from the template fits to this
\bibitem{novakNS71}
Maung~Zho~Newn and I.~F.\ Sharygin.
\newblock Optimal cubature formulas in the classes $D_2^{1,c}$ and 
$D_2^{1,l_1}$.
\newblock In \emph{Problems of Numerical and Applied Mathematics}, pages 22--27.
\newblock Institute of Cybernetics, Uzbek Academy of Sciences, 1991, in Russian. 

\bibitem{novakNUU15}
Van~Kien\ Nguyen, M.~Ullrich, and T.~Ullrich.
\newblock Boundedness of pointwise multiplication and change of variable and 
applications to numerical integration.
\newblock In preparation.

\bibitem{novakNi92}       
H.~Niederreiter.
\newblock \emph{Random Number Generation and Quasi-Monte Carlo Methods}.
\newblock SIAM, 1992.

\bibitem{novakNHMR11}
B.~Niu, F.~Hickernell, Th.~M\"uller-Gronbach, and K.~Ritter.
\newblock Deterministic multi-level algorithms for infinite-dimensional 
integration on $\R^\N$.
\newblock \novakJC 27:331--351, 2011.  

\bibitem{novakNo88}       
E.~Novak.
\newblock \emph{Deterministic and Stochastic Error Bounds in      
Numerical Analysis}.
\newblock Lecture Notes in Mathematics 1349, Springer-Verlag, 1988.        

\bibitem{novakNo96}       
E.~Novak.
\newblock On the power of adaption.
\newblock \novakJC 12:199--237, 1996.        

\bibitem{novakN01}       
E.~Novak.
\newblock Quantum complexity of integration.
\newblock \novakJC 17:2--16, 2001.       

\bibitem{novakNR96b}
E.~Novak and K.~Ritter.
\newblock High dimensional integration of smooth functions over cubes.
\newblock \novakNM 75:79--97, 1996.        

\bibitem{novakNR97}       
E.~Novak and K.~Ritter.
\newblock The curse of dimension and a       
universal method for numerical integration.
\newblock In G.~N\"urnberger, J.~W.\ Schmidt, and G.~Walz,
editors, \emph{Multivariate Approximation and Splines},
pages 177--188.
ISNM {\bf 125}, Birkh\"auser, 1997.

\bibitem{novakNR99}       
E.~Novak and K.~Ritter.
\newblock Simple cubature formulas with high polynomial exactness.
\newblock \novakCA 15:499--522, 1999.        

\bibitem{novakNR14}
E.~Novak and D.~Rudolf.
\newblock Computation of expectations by Markov chain 
Monte Carlo methods.
\newblock In S.~Dahlke et al., editors,
\emph{Extraction of quantifiable information from complex systems}.
Springer-Verlag, 2014.

\bibitem{novakNSW97}       
E.~Novak, I.~H.\ Sloan, and H.~Wo\'zniakowski.
\newblock Tractability of tensor product linear operators.
\newblock \novakJC 13:387--418, 1997.        

\bibitem{novakNT06}       
E.~Novak and  H.~Triebel.
\newblock Function spaces in Lipschitz domains and optimal rates of convergence        
for sampling.
\newblock \novakCA 23:325--350, 2006.        

\bibitem{novakNUW13}
E.~Novak, M.~Ullrich, and H.~Wo\'zniakowski.
\newblock Complexity of oscillatory integration 
for univariate Sobolev spaces.
\newblock \novakJC 31:15--41, 2015. 

\bibitem{novakNW99}    
E.~Novak and H.~Wo\'zniakowski.
\newblock Intractability results for integration and discrepancy.
\newblock \novakJC 17:388--441, 2001.       

\bibitem{novakNW08}       
E.~Novak and H.~Wo\'zniakowski.
\newblock \emph{Tractability of Multivariate Problems},        
Volume I: Linear Information.
\newblock European Mathematical Society, 2008.

\bibitem{novakNW10}      
E.~Novak and H.~Wo\'zniakowski.             
\newblock \emph{Tractability of Multivariate Problems},      
Volume II: Standard Information for Functionals.      
\newblock European Mathematical Society, 2010.       

\bibitem{novakNW11}
E.~Novak and H.~Wo\'zniakowski.
\newblock Lower bounds on the complexity for linear functionals   
in the randomized setting.
\newblock \novakJC 27:1--22, 2011.   

\bibitem{novakNW12}
E.~Novak and H.~Wo\'zniakowski.
\newblock \emph{Tractability of Multivariate Problems},      
Volume III: Standard Information for Operators.
\newblock European Mathematical Society, 2012.       

\bibitem{novakNC06a}
D.~Nuyens and R.~Cools.
\newblock Fast algorithms for component-by-component construction       
of rank-$1$ lattice rules in shift invariant reproducing kernel Hilbert       
spaces.
\newblock \novakMC 75:903--920, 2006.       

\bibitem{novakNC06b}       
D.~Nuyens and R.~Cools.
\newblock Fast algorithms for component-by-component construction       
of rank-$1$ lattice rules with a non-prime number of points.
\newblock \novakJC 22:4--28, 2006.       

\bibitem{novakPW09}
L.~Plaskota and G.~W.\ Wasilkowski.
\newblock The power of adaptive algorithms for functions with singularities.
\newblock \emph{Journal of Fixed Point Theory and Applications} 6:227--248, 2009. 

\bibitem{novakPW11}
L.~Plaskota and G.~W.\ Wasilkowski. 
\newblock Tractability of infinite-dimensional integration in the worst case and 
randomized settings.
\newblock \novakJC 27:505--518, 2011.

\bibitem{novakRo54}
K.~F.\ Roth.
\newblock On irregularities of distributions.
\newblock \emph{Mathematika} 1:73--79, 1954.        

\bibitem{novakRo80}       
K.~F.\ Roth.
\newblock On irregularities of distributions IV.
\newblock \emph{Acta Arithmetica} 37:67--75, 1980.         

\bibitem{novakRu12}
D.~Rudolf.
\newblock Explicit error bounds for Markov chain Monte Carlo.
\newblock \emph{Dissertationes Mathematicae} 485, 2012. 

\bibitem{novakSU06}
W.~Sickel and T.~Ullrich.
\newblock Smolyak's algorithm, sampling on sparse grids and function spaces        
of dominating mixed smoothness.
\newblock \emph{East Journal on Approximation } 13:387--425, 2007.        

\bibitem{novakSU11}
W.~Sickel and T.~Ullrich.
\newblock Spline interpolation on sparse grids.
\newblock \emph{Applicable Analysis} 90:337--383, 2011.

\bibitem{novakSk95}
M.~M.\ Skriganov.
\newblock Constructions of uniform distributions in terms of 
geometry of numbers.
\newblock \emph{St. Petersburg Mathematical Journal } 6:635--664, 1995.  

\bibitem{novakSKJ02}        
I.~H.\ Sloan, F.~Y.\ Kuo, and S.~Joe.
\newblock On the step-by-step construction of quasi-Monte Carlo       
integration rules that achieves strong tractability error       
bounds in weighted Sobolev spaces.
\newblock \novakMC 71:1609--1640, 2002.        
  
\bibitem{novakSR02}
I.~H.\ Sloan and A.~V.\ Reztsov.
\newblock Component-by-component construction of good lattice rules.     
\newblock \novakMC 71:263--273, 2002.

\bibitem{novakSWW}
I.~H.\ Sloan, X.~Wang, and H.~Wo\'zniakowski.
\newblock Finite-order weights imply tractability of multivariate        
integration.
\newblock \novakJC 20:46--74, 2004.       

\bibitem{novakSW98}       
I.~H.\ Sloan and H.~Wo\'zniakowski.
\newblock When are quasi-Monte Carlo algorithms efficient for high       
dimensional integrals?       
\newblock \novakJC 14:1--33, 1998.

\bibitem{novakSm63}        
S.~A.\ Smolyak.
\newblock Quadrature and interpolation formulas for tensor products       
of certain classes of functions.
\newblock \emph{Doklady Akademy Nauk  SSSR} 4:240--243, 1963.        

\bibitem{novakSu79}
A.~G.\ Sukharev.    
\newblock Optimal numerical integration formulas for some classes of functions.
\newblock \emph{Soviet Mathematics Doklady} 20:472--475, 1979.        

\bibitem{novakTem87}       
V.~N.\ Temlyakov.
\newblock Approximate recovery of periodic functions of several       
variables.
\newblock \emph{Mathematics  USSR Sbornik } 56:249--261, 1987.

\bibitem{novakTem90}       
V.~N.\ Temlyakov.
\newblock On a way of obtaining lower estimates for the error of       
quadrature formulas.
\newblock \emph{Math. USSR Sb.} 181:1403-1413, 1990, in Russian.
\newblock English translation: \emph{Mathematics USSR Sbornik} 71:247--257, 1992.       

\bibitem{novakT93}
V.~N.\ Temlyakov.
\newblock On approximate recovery of functions with bounded mixed derivative.  
\newblock \novakJC 9:41--59, 1993.  
  
\bibitem{novakTem03a}        
V.~N.\ Temlyakov.   
\newblock Cubature formulas, discrepancy, and nonlinear approximation.
\newblock \novakJC 19:352--391, 2003.        

\bibitem{novakTWW88}        
J.~F.\ Traub, G.~W.\ Wasilkowski, and H.~Wo\'zniakowski.
\newblock \emph{Information-Based Complexity}.
\newblock Academic Press, 1988.       

\bibitem{novakTW80}        
J.~F.\ Traub and H.~Wo\'zniakowski.
\newblock \emph{A General Theory of Optimal Algorithms}.
\newblock Academic Press, 1980.       

\bibitem{novakTW02}        
J.~F.\ Traub and H.~Wo\'zniakowski.
\newblock Path integration on a quantum computer.       
\newblock \emph{Quantum Information Processing } 1:365--388, 2003.

\bibitem{novakTr10}   
H.~Triebel.
\newblock \emph{Bases in Function Spaces, Sampling, Discrepancy,    
Numerical Integration}.
\newblock European Mathematical Society, 2010.

\bibitem{novakUl14}
M.~Ullrich.
\newblock On ``Upper error bounds for quadrature formulas 
on function classes'' by K. K. Frolov.
\newblock Preprint, 2014.

\bibitem{novakUU15}
M.~Ullrich and T.~Ullrich.
\newblock The role of Frolov's cubature formula for functions 
with bounded mixed derivative.
\newblock Manuscript, available at arXiv. 

\bibitem{novakVy06}       
J.~Vyb{\'\i}ral.
\newblock Sampling numbers and function spaces.
\newblock \novakJC 23:773--792, 2007.        

\bibitem{novakVyb13}
J.~Vyb{\'\i}ral.
\newblock Weak and quasi-polynomial tractability of approximation
of infinitely differentiable functions.
\newblock \novakJC 30:48--55, 2014. 

\bibitem{novakWa14}
G.~W.\ Wasilkowski.
\newblock Average case tractability of approximating 
$\infty$-variate functions.
\newblock \novakMC 83:1319--1336, 2014.

\bibitem{novakWW95}
G.~W.\ Wasilkowski and H.~Wo\'zniakowski.
\newblock Explicit cost bounds of algorithms
for multivariate tensor product problems.
\newblock \novakJC 11:1--56, 1995.        

\bibitem{novakWW96}       
G.~W.\ Wasilkowski and H.~Wo\'zniakowski.
\newblock On tractability of path integration.
\newblock \emph{Journal Mathematical Physics } 37:2071--2088, 1996.

\bibitem{novakWW99}        
G.~W.\ Wasilkowski and H.~Wo\'zniakowski.
\newblock Weighted tensor-product algorithms for linear multivariate problems.
\newblock \novakJC 15:402--447, 1999.

\end{thebibliography}
% and then copy paste the contents of the .bbl file here for the final version.
%
% E.g.:

\end{document}